\documentclass{elsarticle}

\usepackage{graphicx}
\usepackage[margin=2cm,a4paper]{geometry}
\usepackage{amssymb}
\usepackage{amsmath}
\usepackage{multirow,booktabs,subcaption}
\usepackage{todonotes}
\usepackage[english,francais]{babel}
\usepackage[utf8]{inputenc}
\usepackage[T1]{fontenc}

\newtheorem{e-proposition}[theorem]{Proposition}

\newtheorem{e-definition}[theorem]{Definition\rm}
\newtheorem{remark}{\it Remark\/}

\newcommand{\bmu}{\boldsymbol{\mu}}
\newcommand{\R}{\mathbb{R}}
\newcommand{\D}{\ensuremath {\mathcal D} }
\newcommand{\ds}{\displaystyle}
\newcommand{\bx}{\ensuremath {\boldsymbol{x}} }
\newcommand{\bxi}{\ensuremath {\boldsymbol{\xi}} }
\newcommand{\bt}{\ensuremath {\boldsymbol{t}} }
\newcommand{\br}{\ensuremath {\boldsymbol{r}} }
\newcommand{\bq}{\ensuremath {\boldsymbol{q}} }
\newcommand{\N}{ \ensuremath{\mathcal N} }
\newcommand{\eim}[1]{\bar{#1}}

\def\og{\leavevmode\raise.3ex\hbox{$\scriptscriptstyle\langle\!\langle$~}}
\def\fg{\leavevmode\raise.3ex\hbox{~$\!\scriptscriptstyle\,\rangle\!\rangle$}}

\journal{the Acad\'emie des sciences}
\begin{document}
% place in the next line the header (rubrique) chosen for your article,
% if you know it (you can also have 2, format : Header1/Header2
\centerline{}
\begin{frontmatter}


% \author{Name\thanksref{label2}}
\selectlanguage{english}
\title{Simultaneous Empirical Interpolation and Reduced Basis method
  for non-linear problems}

% The [label1] can be suppressed if there is only one address for all authors

\selectlanguage{english}
\author{C\'ecile Daversin},
\ead{daversin@math.unistra.fr}
% \author[mit]{Anthony T. Patera},
% \ead{patera@mit.edu}
\author{Christophe Prud'homme}
\ead{prudhomme@unistra.fr}

\address{IRMA (UMR 7501) 7 rue Ren\'e-Descartes 67084 Strasbourg Cedex}
% \address[mit]{Massachusetts Institute of Technology, Department of Mechanical Engineering, Romm 3-264, 77 Massachusetts Avenue
% Cambridge MA 02139-4307, USA}

% If you know the dates of reception, and acceptation you can put them now;
%  idem the name of the person presenting the Note

\medskip
\begin{center}
{\small Received *****; accepted after revision +++++\\
Presented by £££££}
\end{center}

\begin{abstract}
\selectlanguage{english}

In this paper, we focus on the reduced basis methodology in the
context of non-linear non-affinely parametrized partial differential
equations in which affine decomposition necessary for the reduced
basis methodology are not
obtained~\cite{veroy03:_poster_error_bound_reduc_basis,veroy03:_reduc}.
To deal with this issue, it is now standard to apply the EIM
methodology~\cite{EIM,EIMgrepl} before deploying the Reduced Basis
(RB) methodology. However the computational cost is generally huge as
it requires many finite element solves, hence making it inefficient,
to build the EIM approximation of the non-linear
terms~\cite{EIMgrepl,dptv13}. We propose a simultaneous EIM Reduced
basis algorithm, named SER, that provides a huge computational gain
and requires as little as $N+1$ finite element solves where $N$ is the
dimension of the RB approximation. The paper is organized as follows:
we first review the EIM and RB methodologies applied to non-linear
problems and identify the main issue, then we present SER and some
variants and finally illustrates its performances in a benchmark
proposed in~\cite{EIMgrepl}.  {\it To cite this article: C. Daversin,
  C. Prud'homme, C. R. Acad. Sci. Paris, Ser. I 340 (2015).}

\vskip 0.5\baselineskip

\selectlanguage{francais}

\noindent{\bf R\'esum\'e} \vskip 0.5\baselineskip \noindent
{\bf Une méthode EIM Bases Réduites simultanée pour les équations aux
  dérivées partielle non-linéaires et non-affines. }

Dans ce papier, nous nous intéressons à la méthodologie bases réduites
(RB) dans le contexte d'équations aux dérivées partielles
paramétrisées non-linéaires et non-affines et pour lesquelles la
décomposition affine nécessaire à la méthodologie RB ne peut être
obtenue~\cite{veroy03:_poster_error_bound_reduc_basis,veroy03:_reduc}.
Pour traiter ce problème, il est à présent standard d'appliquer la
méthodologie EIM~\cite{EIM,EIMgrepl} avant de déployer la méthodologie
RB. Cependant le coût de calcul de cette approche est en général
considérable car il requiert de nombreuses évaluations élément fini,
la rendant très peu compétitive, pour construire l'approximation EIM
des termes non-linéaires~\cite{EIMgrepl,dptv13}. Nous proposons
l'algorithme SER qui construit simultanément l'approximation EIM et
RB, fournit ainsi un gain de calcul considérable et requiert au
minimum $N+1$ résolutions élément fini où $N$ est la dimension de
l'approximation RB. Le papier est organisé comme suit: tout d'abord
nous passons en revue les méthodes EIM et RB appliquées aux problèmes
non-linéaires et identifions la difficulté principale, puis nous
présentons SER et quelques variantes et finalement nous illustrons
ses performances sur un benchmark proposé par~\cite{EIMgrepl}.

{\it Pour citer cet article~: C. Daversin, C. Prud'homme, C. R. Acad. Sci. Paris, Ser. I 340 (2015).}

\end{abstract}
\end{frontmatter}

% now the Version française abrégée, if it exists
\selectlanguage{francais}
%\section*{Version fran\c{c}aise abr\'eg\'ee}
% Text of your Version française abrégée here.
% Note you do not need to repeat here equations that you use in the
% main text - for example 'voir (3)' is quite acceptable.

\selectlanguage{english}
% main text
\section{Introduction}
\label{sec:introduction}

%% Reduced Basis Method
Reduced order modeling is more and more used in engineering problems
due to efficient evaluation of quantities of interest.  The Reduced
Basis Method (see
\cite{prud'homme:70,veroy03:_reduc,veroy03:_poster_error_bound_reduc_basis,MR2061274,Reviewquarteroni2011certified,rozza2007reduced})
has been especially designed for real-time and many-query contexts,
and cover a large range of problems among which non-affinely
parametrized Partial Differential Equations (PDE). A core enabler of this
method is the so-called offline/online decomposition of the
problem. This allows for computing costly parameter independent terms
that depend solely on the finite element dimension. However such
decomposition is not necessarily or readily available in particular
for non-affine/non-linear problems. The Empirical Interpolation Method
(EIM, see \cite{EIM,EIMgrepl}) has been developed to recover this core
ingredient and is used prior to the reduced basis methodology on
industrial based applications, see e.g. \cite{dptv13}. However the EIM
building step can be costly when the terms are non-linear and
requires many non-linear finite element solves. It is a deterring
trait for this methodology which forbids its application to non-linear
applications.

In this paper, we propose a Simultaneous EIM-RB (named SER)
construction that requires only but a few finite element solves and
which builds together the affine decomposition as well as the RB
ingredient.  To start we first give an overview of both Reduced Basis
and Empirical Interpolation methods in non-affinely parametrized
PDE context to recall the necessary
notions. Based on these considerations, a second part makes an
assessment of the proposed simultaneous approach discussing the
changes to be made in the EIM offline step. The last part illustrates
our method with preliminar results obtained on a benchmark introduced
in \cite{EIMgrepl}.

\section{Preliminaries} %% RB / EIM description

%% Context
Let $u(\bmu)$ be a solution of a non-linear parametrized PDE, where
$\bmu$ is the $p$-vector of inputs ($\bmu \in \D \subset \R^p$
~\text{the parameter space}). Considering the Hilbert space $X \equiv
H^1_0(\Omega)$, the variational formulation of the PDE consists in
finding $u(\bmu) \in X$ as a root of a functional $r$ such that
\begin{equation}
  \label{eq:1}
  r( u(\bmu),v; \bmu ) = 0 ~\forall v \in X.
\end{equation}
We denote $\N$ the dimension of the finite element approximation space
$X^{\N} \subset X$, and $u_{\N}(\bmu)$ the associated finite element
approximation of $u(\bmu)$ solution of \eqref{eq:1}.

In the following, we use a Newton algorithm to deal with the
non-linearity. Denoting ${}^ku(\bmu)$ the solution at $k^{th}$
Newton's iteration, and $j$ the Jacobian associated with the functional
$r$, the problem consists in finding ${\delta}^{k+1}u(\bmu)\equiv
{}^{k+1}u(\bmu) - {}^ku(\bmu) \in X^\N$ and thus ${}^{k+1}u(\bmu) \in X^\N$ such that
\begin{equation}
\label{equ:newton}
j( {\delta}^{k+1}u(\bmu),v; \bmu ; {}^ku(\bmu) )  = -r( {}^ku(\bmu),v;
\bmu ) ~\forall v \in X^\N
\end{equation}
at each iteration $k$.

%If it exists, the affine decomposition of bilinear and linear forms $a$ and $f$ then reads :
% \begin{equation}
%   \label{equ:affine-dec}
%   \ds{ a(u, v; \bmu) = \sum \limits_{q=1}^{Q_a} \theta^q_a(\bmu) a^q(u, v) } \quad
%   \ds{f(v; \bmu) = \sum \limits_{q=1}^{Q_f} \theta^q_f(\bmu) f^q(v) } \quad \forall u, v \in X, ~\forall \bmu \in \D
% \end{equation}
If $j$ and $r$ depend affinely in $\bmu$, there exists positive integers $Q_j$ and $Q_r$, $\bmu$-dependent functions
$\theta^{q_j}_j$, $\theta^{q_r}_r$ and functions $j^{q_j}$, $r^{q_r}$ ($1 \leqslant q_r \leqslant Q_r$, $1 \leqslant q_j \leqslant Q_j$)
allowing to write such affine decompositions
%the affine decomposition of the Jacobian $j$ and of the residual function $r$ then read :
\begin{equation}
  \label{equ:affine-dec}
  \ds{ j(u, v; \bmu; {}^ku) = \sum \limits_{q_j=1}^{Q_j} \theta^{q_j}_j({}^ku; \bmu)  j^{q_j}(u, v) } \quad
  ~\ds{ r({}^ku, v; \bmu) = \sum \limits_{q_r=1}^{Q_r} \theta^{q_r}_r({}^ku; \bmu) r^{q_r}(v) } \quad \forall u, v \in X^\N, ~\forall \bmu \in \D
\end{equation}

For non-affinely parametrized problems, $j$ and $r$ depend on a set
$\{ w_i(u(\bmu),\bx;\bmu)\}_i$ of non-affine functions
\begin{equation}
  \label{equ:var-with-w}
  %a( u(\bmu),v; \bmu; \{ w_i(u,\bx;\bmu)\}_i) = f(v; \bmu; \{ w_i(u,\bx;\bmu)\}_i) ~\forall v \in X
  j( u, v; \bmu ; {}^ku(\bmu); \{ w_i(u(\bmu),\bx;\bmu)\}_i ) {\delta}^{k+1}u(\bmu) = -r( {}^ku(\bmu),v; \bmu; \{ w_i(u(\bmu),\bx;\bmu)\}_i )
\end{equation}
and the affine decomposition (\ref{equ:affine-dec}) doesn't exist owing to $w_i(u(\bmu),\bx;\bmu)$ definition. In order to recover an affine decomposition for (\ref{equ:var-with-w}), we build an affine approximation $w_{i,M_i}(u(\bmu),\bx;\bmu)$ of $w_i(u(\bmu),\bx;\bmu)$
\begin{equation}
  w_{i,M_i}(u(\bmu),\bx;\bmu) = \sum \limits_{m=1}^{M_i} \beta^{M_i}_{i,m}(u(\bmu);\bmu) q_{i,m}(\bx) ~\forall i
\end{equation}

%% EIM
\subsection{Empirical Interpolation Method} \label{subsec:eim-method}

% \textit{Note : For clarity reasons, EIM ingredients are distinguished by a bar $\eim{(.)}$ in order to avoid confusions with Reduced Basis Method notations in next sections. \\}

In order to build an affine approximation $w_{M}(u(\bmu),\bx;\bmu) = \sum
\limits_{m=1}^{M} \beta^{M}_{m}(u(\bmu);\bmu) q_{m}(\bx)$ of a non-affine
parameter dependent function $w$, we first introduce a sample
$\eim{S}_M = \{ \eim{\bmu}_1, \dots, \eim{\bmu}_M \} \in \D^M$, and
the associated fonction space $\eim{W}_M = span\{ \eim{\bxi}_m \equiv
w(u(\eim{\bmu}_m),\bx;\eim{\bmu}_m), 1 \leqslant m \leqslant M \}$.
These sets are built from a subset $\Xi$ of $\D$, in which the first
sample point $\eim{\bmu}_1$ is picked (assuming that $\eim{\bxi}_1
\neq 0$).  The basis functions $q_m$ are based on $\eim{\bxi}_m$, and
establishing of $\beta^M_m$ coefficients requires, online, the
solution of a $M\times M$ system, ensuring that the EIM approximation
$w_M(u,\bx;\bmu)$ is exact on a set of interpolation points $\{ t_i
\}_i$.
\begin{equation}
  \label{equ:eim-init-offline}
  \eim{S}_1 = \{ \eim{\bmu}_1 \}, \quad
  \eim{\bxi}_1 = w(u(\eim{\bmu}_1),\bx;\eim{\bmu}_1), \quad
  \eim{W}_1 = span\{ \eim{\bxi}_1 \} \quad
\bt_1 = arg \sup_{\bx \in \Omega} |\eim{\bxi}_1(\bx)|, \quad
q_1 = \frac{\eim{\bxi}_1(\bx)}{\eim{\bxi}_1(\bt_1)}
\end{equation}

For $M \geqslant 2$, the sample points $\eim{\bmu}_M$ are determined from a Greedy algorithm :
\begin{equation}
  \label{equ:eim-best-fit-fem}
  \begin{array}{ll}
    \eim{\bmu}_M = arg \max_{\bmu \in \Xi} ~ \inf_{z \in W_{M-1}} ||w(.;.;\bmu) - z ||_{L^{\infty}(\Omega)} \\
    \eim{\bxi}_M = w \left( u(\eim{\bmu}_M), \bx ; \eim{\bmu}_M  \right),
    \quad \eim{S}_M = \eim{S}_{M-1} \cup \{ \eim{\bmu}_M \},
    \quad \eim{W}_M = \eim{W}_{M-1} \oplus span \{ \eim{\bxi}_M \}
  \end{array}
\end{equation}

The system ensuring the exactness of $w_{M-1}$ at
$\{t_i\}_{i=1}^{M-1}$, gives
$w_{M-1}(u(\eim{\bmu}_M),\bx;\eim{\bmu}_M)$, leading to the residual
$\br_M$ defined as $\br_M(\bx) = w(u(\eim{\bmu}_M),\bx;\eim{\bmu}_M) -
w_{M-1}(u(\eim{\bmu}_M),\bx;\eim{\bmu}_M)$ on which interpolation
points $\{ \bt_M \}_{M \leqslant 2}$ and basis functions $\{ q_M \}_{M
  \leqslant 2}$ are computed. The next $\bt_M$ and $\bq_M$ are given by
\begin{equation}
  \bt_M = arg \sup_{\bx \in \Omega} |\br_M(\bx)|, \quad
  \bq_M(\bx) = \frac{\br_M(\bx)}{\br_M(\bt_M)}
\end{equation}

We are now ready to apply EIM on the jacobian $j$ the affine
decomposition reads
\begin{equation}
  \label{equ:affine-dec-w}
    % \ds{ a(u, v; \bmu; \{ w_i(u(\bmu),\bx;\bmu) \}_i) =
    %   \underbrace{\sum \limits_{q=1}^{Q_a^{eim}} \sum \limits_{m=1}^{M_q^a} \beta^q_{a,m}(u,\bmu) a_{m}^{q}(u,v) }_\textrm{non-affine part of $a$}
    %   + \underbrace{\sum \limits_{l=1}^{Q_a^{aff}} \theta^l_a(\bmu) a^l(u, v)}_\textrm{affine part of $a$}
    \ds{ j(u, v; \bmu; {}^ku; \{ w_i(u(\bmu),\bx;\bmu) \}_i) =
      \underbrace{\sum \limits_{q=1}^{Q_j^{eim}} \sum \limits_{m=1}^{M_q^j} \gamma^q_{j,m}({}^ku; \bmu) j_{m}^{q}(u,v) }_\textrm{non-affine part of $j$}
      + \underbrace{\sum \limits_{l=1}^{Q_j^{aff}} \theta^l_j({}^ku; \bmu) j^l(u, v)}_\textrm{affine part of $j$}
     \quad \forall u, v \in X, ~\forall \bmu \in \D }, \\
\end{equation}

with $Q_j^{eim}$ (resp. $Q_j^{aff}$) the number of non-affinely
(resp. affinely) parametrized terms of $j$,
and $\gamma_{j,m}^q$ (resp. $j_{m}^{q}$) obtained from
$\beta^{M_i}_{i,m}$ (resp. $q_{i,m}$).
Similar decomposition is obtained for the residual $r$.

%% RBM
\subsection{Reduced basis method} \label{subsec:rb-method}

We now turn to the RB methodology. We introduce $S_N = \{\bmu_1,
\cdots, \bmu_N\}$ with $N << \N$, and we define the
set of solutions $S_N^u = \{ u_{\N}(\bmu_i) \}_{i=1}^N$ which are
orthonormalized with respect to the $<,>_X$ inner
product to provide $W_N = \{ \bxi_n
\}_{i=1}^N$.

The reduced basis approximation $u_N(\bmu) \in W_N$ is expressed as a linear
combination of $W_N$ elements
\begin{equation}
  \label{equ:reduc-approx-def}
 %  u_N(\bmu) = \sum_{i=1}^N u_{N,i}(\bmu) \bxi_i ~\mid
 % ~\ds{ a(u_{N}, v; \bmu; \{ w_i(u_{N},\bx;\bmu) \}_i) = f(v; \bmu; \{ w_i(u_{N},\bx;\bmu) \}_i) ~\forall v \in X^{\N} }
  u_N(\bmu) = \sum_{i=1}^N u_{N,i}(\bmu) \bxi_i ~\mbox{ such that }
 ~\ds{ r(u_{N}, v; \bmu; \{ w_i(u_{N},\bx;\bmu) \}_i) = 0 ~\forall v \in W_N}
\end{equation}

The offline/online strategy used to build efficiently the reduced
basis approximation is based on the affine decomposition
(\ref{equ:affine-dec}),(\ref{equ:affine-dec-w}). Choosing $\{ \bxi_n
\}_{i=1}^N$ as test functions and using (\ref{equ:reduc-approx-def})
we have for $1 \leqslant l \leqslant N$
% \begin{equation}
%   \ds{ \sum \limits_{j=1}^N a(\bxi_j,\bxi_k; \bmu; \{ w_i(u(\bmu),\bx;\bmu) \}_i) U_{N,j} = f(\bxi_k; \bmu; \{ w_i(u(\bmu),\bx;\bmu) \}_i), \quad 1 \leqslant k \leqslant N}
% \end{equation}
\begin{equation}
\label{eq:2}
  \ds{ \sum \limits_{j=1}^N j(\bxi_j,\bxi_l; \bmu; {}^ku_N; \{ w_i({}^ku_N,\bx;\bmu) \}_i)
    ({}^{k+1}u_{N,j} - {}^{k}u_{N,j}) =
    -r(\bxi_l; \bmu; {}^ku_N; \{ w_i({}^ku_N,\bx;\bmu) \}_i)}
\end{equation}

% which gives, using the affine decomposition (\ref{equ:affine-dec-w}) with $ \leqslant k \leqslant N $
% \begin{equation*}
%   %\label{equ:rb-online-system-w}
%   \sum \limits_{j=1}^N \left[
%     \sum \limits_{q=1}^{Q_a^{eim}} \sum \limits_{m=1}^{M_q^a} \beta^q_{a,m}(u,\bmu) a_{m}^q(\bxi_j, \bxi_k)
%     +\sum \limits_{l=1}^{Q_a^{aff}} \theta^l_a(\bmu) a^l(\bxi_j, \bxi_k)
%   \right] u_{N,j} =
%   \sum \limits_{q=1}^{Q_f^{eim}} \sum \limits_{m=1}^{M_q^f} \beta^q_{f,m}(u,\bmu) f_{m}^q(\bxi_k)
%   +\sum \limits_{l=1}^{Q_f^{aff}} \theta^l_f(\bmu) f^l(\bxi_k)
% \end{equation*}

% \begin{equation*}
% \sum \limits_{j=1}^N \left[
%   \sum \limits_{q=1}^{Q_j^{eim}} \sum \limits_{m=1}^{M_q^j} \beta^q_{j,m}(\bmu, {}^ku) j_{m}^q(\bxi_j, \bxi_k)
%   +\sum \limits_{l=1}^{Q_j^{aff}} \theta^l_j(\bmu, {}^ku) j^l(\bxi_j, \bxi_k)
% \right] \delta^k u_{N,j} =
% - \sum \limits_{q=1}^{Q_r^{eim}} \sum \limits_{m=1}^{M_q^r} \beta^q_{r,m}(\bmu, {}^ku) r_{m}^q(\bxi_k)
% + \sum \limits_{l=1}^{Q_r^{aff}} \theta^l_r(\bmu, {}^ku) r^l(\bxi_k)
% \end{equation*}

The offline step consists in the pre-computation of the terms
$j_{m}^q(\bxi_j, \bxi_l), j^l(\bxi_j, \bxi_l)$ and $r_{m}^q(\bxi_l),
r^l(\bxi_l)$, and is done offline once thanks to the affine
decomposition. To evaluate online ${}^{k+1}u_N(\bmu)$ for any given $\bmu$,
the coefficients $\gamma^q_{j,m}({}^ku_{N},\bmu),
\theta^l_j({}^ku_{N},\bmu)$ and $\gamma^q_{r,m}({}^ku_{N},\bmu),
\theta^l_r({}^ku_{N},\bmu)$ are computed in order to recover the $N
\times N$ system (\ref{equ:reduc-approx-def}) giving ${}^{k+1}u_{N,j}$
coefficients and consequently ${}^{k+1}u_N(\bmu)$.

% As previously mentionned, the EIM offline step is currently performed up front, independently from RB offline step. The following section makes a brief assessment of this approach, and proposes another one combining both offline steps.

\section{A Simultaneous EIM-RB method} %% optimization of the method
\label{sec:co-build-method}

The standard methodology described previously briefly requires
the use of EIM prior to the RB methodology for each affine parameter
dependent function $w_i$ (\ref{equ:var-with-w}). These EIM
approximations are built up front allowing to then write the affine
decomposition~\eqref{equ:affine-dec-w}.  However the greedy algorithm
(\ref{equ:eim-best-fit-fem}) of the EIM offline step requires the
computation of the solution $u_{\N}(\bmu)$ for all points of $\Xi
\in \D$, \emph{i.e} in the step \eqref{equ:eim-best-fit-fem}.

% This solution is approximated by the
% parametric finite element solution $u^p_{\N}(\bmu)$ --based on affine
% decomposition -- and the number of finite element computations is
% consequently proportional to the size of the trainset $\Xi$.  Since
% the trainset $\Xi$ doesn't vary with $M$, the set of solution $\{
% u^p_{\N}(\bmu) \}_{\bmu \in \Xi}$ computed for $M=2$ (the first step
% $M=1$ needs only $u^p_{\N}(\eim{\bmu}_1)$) can be stored and re-used
% for $M \geqslant 3$. The number of finite element solve is then
% actually reduced to the size of $\Xi$ which can though be large and
% leads to a very costly offline step since $\N$ is typically large.

The SER methodology proposes to reduce the computational cost by
simply using the readily available reduced basis approximation based
on the previous EIM step and build simultaneously EIM and RB.  In
particular the expensive step (\ref{equ:eim-best-fit-fem}) will now
use solely RB approximations.

% To wit
% (\ref{equ:eim-best-fit-fem}) using $u_{N}(\bmu)$ now reads at each
% Newton's iteration $k$,
% \begin{equation}
%   \label{equ:best-fit}
%   \eim{\bmu}_{i,M} = arg \max_{\bmu \in \Xi} ~ \inf_{z \in W_{M-1}} ||w_i( {}^ku_{N}(\bmu),.;\bmu) - z ||_{L^{\infty}(\Omega)}
% \end{equation}

% \begin{remark}
% \label{rem:compromise-fem-for-accuracy}
% \textit{The use of ${}^ku_{N}$ can spoil the choice of $\eim{\bmu}_{i,M}$ due to a lack of precision, especially for first basis obtained from a rough EIM approximation. A compromise should be fine to ensure a sufficient precision depending on the complexity of the considered model.\\}
% \end{remark}

% The aim is to build EIM and RB basis concurently. This section
% describes the alternative build of these basis one by one.  This
% process can clearly be generalized to a build per group of
% customizable size.

The initialization of EIM offline stage (\ref{equ:eim-init-offline})
doesn't change since no reduced approximation $u_{N}(\bmu)$ is available
yet.  The resulting rough EIM approximations $w_{i,1}$ are used to
compute a first affine decomposition from which $\bxi_1$ can be built,
leading to a first reduced approximation $u_1(\bmu)$. After this
initialization step, EIM and RB approximation spaces are enriched
alternatively. Each new EIM basis function $q_M$ : \textit{(i)} is
built from reduced basis approximation $u_{M-1}(\bmu)$ obtained at
previous iteration $M-1$ \textit{(ii)} complete the EIM approximation
to then build the affine decomposition for $u_M(\bmu)$.  To
summarize, the EIM steps modified by SER now read
\begin{equation}
  \eim{\bmu}_{i,M} = arg \max_{\bmu \in \Xi} ~ \inf_{z \in W_{M-1}} ||w_i( u_{M-1}(\bmu);.;\bmu) - z ||_{L^{\infty}(\Omega)}
  \label{eq:bestfit-opt}
\end{equation}
\begin{equation}
  \label{eq:eim-basis-computation}
  \eim{\bxi}_{i,M} = w \left( u_{M-1}(\eim{\bmu}_{i,M}); \bx; \eim{\bmu}_{i,M} \right),
  \quad \eim{S}_{i,M} = \eim{S}_{i,M-1} \cup \{ \eim{\bmu}_{i,M} \},
  \quad \eim{W}_{i,M} = \eim{W}_{i,M-1} \oplus span \{ \eim{\bxi}_{i,M} \}
\end{equation}
\noindent
The residual $\br_{i,M}(\bx) = w_i(u_{M-1}(\eim{\bmu}_{i,M}),\bx;\eim{\bmu}_{i,M}) - w_{i,M-1}(u_{M-1}(\eim{\bmu}_{i,M}),\bx;\eim{\bmu}_{i,M})$
then gives $\bt_{i,M}$ and $q_{i,M}$
\begin{equation}
  \bt_{i,M} = arg \sup_{\bx \in \Omega} |\br_{i,M}(\bx)|, \quad
  \bq_{i,M}(\bx) = \frac{\br_{i,M}(\bx)}{\br_{i,M}(\bt_{i,M})}
\end{equation}

% \begin{remark}
% \label{rem:use-fem-in-residual}
% \textit{The reduced basis approximation $u_{M-1}$ is used in evaluation of $w_i$ to compute $\eim{\bxi}_{i,M}$ (\ref{eq:eim-basis-computation}) and in residual computation. The gain in term of computational time shouldn't be significant since this computation already depend on space and consequently on finite element dimension $\N$. This computation being done once, use of $u^p_{\N}$ could ensure more accuracy without being too costly.} \\
% \end{remark}

The affine decomposition is then updated with the new EIM
approximation $w_{i,M}$, to compute the next reduced basis $\bxi_M$
(\emph{i.e}) thanks to \eqref{eq:2} to then enrich $W_N$ space.\\

% to find the parametric finite element solution $u_{\N}^p$ of the system
% \begin{equation}
% \label{equ:offline-fixed-point}
% \left( \sum \limits_{m=1}^M \beta_{q,m}^a(u_{\N}^p(\bmu_M), \bmu_M) a_{q,m} \right) u_{\N}^p(\bmu_M)
% = \sum \limits_{m=1}^M \beta_{q,m}^f(u_{\N}^p(\bmu_M), \bmu_M) f_{q,m}
% \end{equation}
% leading to $\bxi_M$ to enrich $W_M$ space.\\

\begin{remark}[RB updates]
  \label{remark:rebuild-rb}
  The EIM approximations are changing during the build of
    $W_N$ space. Then, the operator -- and the solved problem --
    evolves at each step, which can mildly deteriorate the approximation.
    %as seen in the numerical experiments.
    Recomputing $W_N$ elements
    using finite element solves based on the current EIM approximation
    may be considered. The offline precomputations need however to be
    updated for all elements of $W_N$ for each update of the affine
    decomposition.\\
\end{remark}
\begin{remark}[EIM updates]
  \label{remark:simultaneous-eim-rb}
  In the SER methodology only one finite element solve of \eqref{eq:1}
  is required for the initialization afterwards the EIM and RB
  approximations are updated alternatively. We may in fact update the
  RB approximation every $r$ EIM steps: if $r=M$ we recover the
  standard method, if $r=1$ we recover the SER method and if $1 < r <
  M$ we have an intermediairy method which requires finite element
  solves for the intensive EIM step \eqref{equ:eim-best-fit-fem}
  before the first RB update. Other alternatives are readily available
  and some will be discussed in a future publication.
\end{remark}

\section{Preliminary results} %% first results
\label{sec:results}
%\todo[inline]{ space discretisation ? idem for eim }
We now turn to numerical experiments of the SER method compared to the
standard one based on the benchmark problem introduced in
\cite{EIMgrepl}. It reads, find $u$ such that
\begin{equation}
  \label{equ:model-bench}
  -\Delta u + \mu_1 \frac{e^{\mu_2 u} - 1}{\mu_2} = 100 \sin(2\pi x) \sin(2\pi y)
  ~\text{in} ~\Omega = ]0,1[^2 ~\text{and} ~\bmu=(\mu_1,\mu_2) \in \D = [0.01, 10]^2
\end{equation}
and we are interested in the output $s$, the average of the solution
$u$ over $\Omega$. Thanks to the function $g(u,x;\bmu) = \mu_1
\frac{e^{\mu_2 u} - 1}{\mu_2}$ we are in the setting of the SER
methodology for non-linear problems. Following the standard
methodology we develop $g_M = \sum \limits_{i=1}^M \beta_m^M(u,\bmu)
q_m(\bx)$ of $g$ thanks to the Empirical Interpolation Method (see
section \ref{subsec:eim-method}). We shall use the absolute errors on
the solution and the output, defined as follows
\begin{equation}
  \label{equ:errors}
  \epsilon^{u,r}_{M,N} = \parallel u_{\N} - u^r_N \parallel_{L_2}
  \quad
  \epsilon^{s,r}_{M,N} = \mid s_{\N} - s^r_N \mid
\end{equation}
where $\cdot_{\N}$ is the finite element solution/output of the
initial problem \eqref{equ:newton}, $\cdot^r_N$ the reduced basis
solution/output while $r$ is the frequency at which the EIM are
updated, see remark~\ref{remark:simultaneous-eim-rb}.

We first display the maximum of the absolute errors (\ref{equ:errors})
using the standard method $r=M$ in table \ref{tab:crb-P3} which
reproduces the results in \cite{EIMgrepl}.
% Next we tabulate in tables
% \ref{tab:crb-cobuild-P3-1}, \ref{tab:crb-cobuild-P3-1-no-restart}
% $\epsilon_{M,N}^{u,r}$ and $\epsilon_{M,N}^{s,r}$ for the same values of $(N,M)$
% for the case $r=1$ which corresponds to the SER methodology.
Tables \ref{tab:crb-cobuild-P3-1} and \ref{tab:crb-cobuild-P3-1-no-restart} correspond
to the SER methodology ($r=1$) with $N=M$ and
%Tables \ref{tab:crb-cobuild-P3-1} and \ref{tab:crb-cobuild-P3-1-no-restart}
investigate the influency of $W_N$ recomputation in SER method, see remark \ref{remark:rebuild-rb}.
Table \ref{tab:crb-cobuild-P3-5} is an intermediary stage such that $1 < r < M$ where the EIM basis
functions are built by groups of $r$ ($r=5$). In this last case, since we do not
have access yet to the RB, the first group of $r$ EIM basis is built
from finite element approximations and reduced basis approximations
are used afterwards.  This explains, for example, that $(N=4,M=5)$ is
similar to the first one in table \ref{tab:crb-P3}.

% As to the SER
% case($r=1$), only one one finite element approximation
% $u(\eim{\bmu_1})$ to build the first EIM basis function
% (\ref{equ:eim-init-offline}) and only reduced basis approximations are
% used afterwards.

The maximum of the errors observed with SER method ($r=1$) is slightly higher
than with the case $r=5$ which itself displays resultats slightly
higher than the standard method $r=M$. This behavior is expected. The
results show the pertinence of the SER method and that we can expect
good results within a reasonnable computational budget, since the number
of finite element approximations can then be reduced to $N+1$ (table \ref{tab:crb-cobuild-P3-1-no-restart}).

\begin{table}[htbp]
  \centering
  \begin{subtable}[t]{.4\linewidth}
      \begin{tabular}{cccc}
    \toprule
    $N$ & $M$ & $\max(\epsilon^{u,M}_{M,N})$ & $\max(\epsilon^{s,M}_{M,N})$    \\
    \midrule
    4 & 5 & 7.38e-3 & 5.75e-3 \\
    8 & 10 & 1.01e-3 & 2.34e-4 \\
    12 & 15 & 1.49e-4 & 3.09e-5 \\
    16 & 20 & 2.21e-5 & 1.25e-5 \\
    20 & 25 & 5.88e-6 & 2.82e-6\\
    \bottomrule
  \end{tabular}
  \caption{$r=M$}
  \label{tab:crb-P3}
  \end{subtable}
  \begin{subtable}[t]{0.4\linewidth}
    \centering
    \begin{tabular}{cccc}
      \toprule
      $N$ & $M$ & $\max(\epsilon^{u,5}_{M,N})$ & $\max(\epsilon^{s,5}_{M,N})$ \\
      \midrule
      4 & 5 & 8.21e-3 & 6.31e-3 \\
      8 & 10 & 4.48e-3 & 6.18e-3 \\
      12 & 15 & 2.69e-4 & 2.36e-4 \\
      16 & 20 & 1.48e-4 & 9.31e-5 \\
      20 & 25 & 2.60e-5 & 1.46e-5 \\
      \bottomrule
    \end{tabular}
    \caption{$r=5$}
    \label{tab:crb-cobuild-P3-5}
  \end{subtable}\\

  \begin{subtable}[t]{0.4\linewidth}
    \centering
    \begin{tabular}{cccc}
      \toprule
      $N$ & $M$ & $\max(\epsilon^{u,1}_{M,N})$ & $\max(\epsilon^{s,1}_{M,N})$ \\
      \midrule
      5 & 5 & 9.98e-3 & 7.77e-3 \\
      10 & 10 & 2.32e-3 & 1.86e-3 \\
      15 & 15 & 4.61e-4 & 3.75e-4 \\
      20 & 20 & 2.48e-4 & 2.02e-4 \\
      25 & 25 & 3.51e-5 & 2.33e-5 \\
      \bottomrule
    \end{tabular}
    \caption{$r=1$ \\($W_N$ recomputed)}
    \label{tab:crb-cobuild-P3-1}
  \end{subtable}
  %\hfill
  \begin{subtable}[t]{0.4\linewidth}
    \centering
    \begin{tabular}{cccc}
      \toprule
      $N$ & $M$ & $\max(\epsilon^{u,1}_{M,N})$ & $\max(\epsilon^{s,1}_{M,N})$ \\
      \midrule
      5 & 5 & 1.30e-2 & 1.02e-2 \\
      10 & 10 & 2.20e-3 & 1.50e-3 \\
      15 & 15 & 4.83e-4 & 4.05e-4 \\
      20 & 20 & 2.42e-4 & 1.98e-4 \\
      25 & 25 & 1.50e-5 & 1.24e-5 \\
      \bottomrule
    \end{tabular}
    \caption{$r=1$ \\($W_N$ not recomputed)}
    \label{tab:crb-cobuild-P3-1-no-restart}
  \end{subtable}
  \caption{\textbf{SER algorithm:} Maximum absolute errors on solution $u$ and on output $s$}
\end{table}

% \appendix

% \section{}
% \label{}
\section*{Conclusion}
We have now an algorithm and variants that allow an efficient use of
EIM and RB in the context of non-linear non-affine partial
differential equations. We have already deployed SER in real
applications and the initial results are very promising showing the
same behavior as the benchmark problem of this paper. This will be
reported as well as an analysis of SER  in a subsequent paper with
applications to non-linear multiphysic problems.

\section*{Acknowledgements}
The authors would like to thank A.T. Patera, for the discusion that
initiated the SER algorithm and the subsequent ones, and S. Veys.
They are also thankful for the financial support of the ANR CHORUS and
the LABEX IRMIA(Strasbourg).

\end{document}